\documentstyle[12pt]{article}
\newcommand{\qdet}{{\cal D}_q}
\newcommand{\lmatr}{\left(\begin{array}}
\newcommand{\rmatr}{\end{array}\right)}
\newcommand{\lpar}{\stackrel{\leftarrow}{\partial}}
\newcommand{\rpar}{\stackrel{\rightarrow}{\partial}}

\newcommand{\rdel}{\stackrel{\rightarrow}{\delta_R}}
\newcommand{\anabla}{\stackrel{\gets}{\nabla}}
\newcommand{\hanabla}{\hat{\stackrel{\gets}{\nabla}}}
\newcommand{\rnabla}{\stackrel{\to}{\nabla}}
\newcommand{\hrnabla}{\hat{\stackrel{\to}{\nabla}}}
\newcommand{\bara}{\stackrel{}{\bar{A}}}
\newcommand{\ttheta}{\tilde{\theta}}
\newcommand{\qtr}{\mbox{Tr}_q}

\newcommand{\bom}{\bar{\omega}}
\newcounter{nnn}

\begin{document}

\begin{center}
{\Large {\bf Matched differential calculus on the quantum groups $
GL_q(2,C),SL_q(2,C),C_q(2|0)$.}}

Akulov V.P., Gershun V.D.

{\sl Kharkov Institute of Physics \& Technology}
\end{center}

\begin{quotation}
We proposed the construction of the differential calculus on the quantum
group and its subgroup with the property of the natural reduction: the
differential calculus on the quantum group  $GL_q(2,C)$ has to contain the
differential calculus on the quantum subgroup $SL_q(2,C)$ and
quantum plane $C_q(2|0)$ (''quantum matrjoshka''). We found, that there are
two differential calculi, associated to the left differential Maurer--Cartan
1-forms and to the right differential 1-forms. Matched reduction take the
degeneracy between the left and right differentials. The classical limit
($q\to 1$) of the ''left'' differential calculus and of the ''right''
differential calculus is undeformed differential calculus.  The condition
${\cal D}_qG=1$ gives the differential calculus on $SL_q(2,C)$, which
contains the differential calculus on the quantum plane $C_q(2|0)$.
\end{quotation}

\section{Introduction}

In this paper we discuss a new effect appearing in the differential calculus
on quantum groups and its subgroups.

After the discovery [1,2] $q$-deformed or quantum groups a peculiar
attention was paid to constructing bicovariant differential calculus on these
quantum groups [3].But bicovariant differential calculus on $GL_q(N)$ cannot
be expressed in terms of differential calculus on $SL_q(N)$, because the
Maurer--Cartan 1-forms have $q^2$-commutation relations with quantum
determinant ${\cal D}_q$. Such an approach was used for quantum groups in
[4-8] with undeformed classical Leibnitz rule for the exterior derivative.
Another approach can be formulated following the Faddeev--Pjatov idea that
the exterior derivative obeys the modified version of the Leibnitz rule [9].

We are going to describe another approach with main assumptions: 1)
undeformed classical Leibnitz rule, 2) the differential calculus on the
quantum group $GL_q(2,C)$ has to contain the differential calculus on the
quantum subgroup $SL_q(2,C)$. As the result, we will find that there are two
differential calculus on the quantum group $GL_q(2,C)$, associated to the
left differential forms and to the right differential forms. The classical
limit $(q\rightarrow 1)$ of the ''left'' dif\-fe\-ren\-tial calculus and the
''right'' dif\-fe\-ren\-tial calculus is the undeformed differential
calculus.

The condition ${\det }_qG=1$ gives the differential calculus on $SL_q(2,C)$. If
the parameter $b$ or $c$ is equal to zero, we will find the differential
calculus on the quantum plane $C_q(2|0)$ [6].

Let us briefly discuss the content of the paper. In the second section the
basic notations of the differential calculus on the quantum group are
introduced. In the third section we   will be dealing with the left
1-forms $ \theta $ and the left differential $\delta_L$. We find the
commutation relations for the left 1-forms and the parameters and we
describe the quantum algebras for the vector fields for $GL_q(2,C)$, but
only one of them has the Drinfeld--Jimbo form. This choice fixes the form
of the quantum trace $ \mbox{Tr}_q\theta $. We derive formulas
corresponding to the commutation relations between the parameters and its
differentials for $GL_q(2,C)$, which contains $SL_q(2,C)$ case and
$C_q(2|0)$ case. In section 4 we considered the $SL_q(2,C)$ case which
also containts $C_q(2|0)$ case. In section 5 we introduce the right
1-forms $\omega $ and right differential $\stackrel{\rightarrow}{\delta
_R}$ like in the previous sections and recovered the results [10].

\section{ Differential calculus on quantum group $GL_q(2,C)$.}

We first recall some basic notations about differential calculus on quantum
group. The starting point for our consideration is the Hopf algebras $
Fun\left(GL_q(2,C)\right)$.

Comultiplication, counit and antipode is determined by

$$
\Delta \left( T_k^i\right) =T_l^i\otimes T_k^l
\eqno(2.1)
$$
$$
e\left( T_k^i\right) =\delta _k^i\;\;\;\;\;\;\;
\eqno(2.2)
$$
$$
S\left( T_k^i\right) =\left( T^{-1}\right) _k^i
\eqno(2.3)
$$

Non-commuting matrix entries $T=\left(
\begin{array}{cc}
T_1^1 & T_2^1 \\
T_1^1 & T_2^2
\end{array}
\right) $ satisfies the RTT-relation
$$
R_{kl}^{ij}T_m^kT_n^l=T_l^jT_k^iR_{mn}^{kl}
\eqno(2.4)
$$
and the quantum Yang-Baxter equation
$$
R_{i_2j_2}^{i_1j_1}R_{i_3k_2}^{i_2k_1}R_{j_3k_3}^{j_2k_2}=
R_{j_2k_2}^{j_1k_1}R_{i_2k_3}^{i_1k_2}R_{i_3j_3}^{i_2j_2}
\eqno(2.5)
$$
where
$$
R_{kl}^{ij}=\left(
\begin{array}{cccc}
R_{11}^{11} & R_{12}^{11} & R_{21}^{11} & R_{22}^{11} \\
R_{11}^{12} & R_{12}^{12} & R_{21}^{12} & R_{22}^{12} \\
R_{11}^{21} & R_{12}^{21} & R_{21}^{21} & R_{22}^{21} \\
R_{11}^{22} & R_{12}^{22} & R_{21}^{22} & R_{22}^{22}
\end{array}
\right) =\left(
\begin{array}{cccc}
q & 0 & 0 & 0 \\
0 & 1 & 0 & 0 \\
0 & \lambda & 1 & 0 \\
0 & 0 & 0 & q
\end{array}
\right)
\eqno(2.6)
$$
where the rows and columns are numerated in the order 11, 12, 21, 22, or
equivalently as composite indices 1,2,3,4 and $\lambda =q-q^{-1}$.

An important element is the quantum determinant
$$
{\cal D}_q={\det }_qT=\sum_{i,k=1}^2(-q)^{I(ik)}T_i^1T_k^2
\eqno(2.7)
$$
where $I(ik)$ is the number of inversions, or through the $q$-deformed
Levi-Civita tensor
$$
\epsilon _q^{ik}=\left(
\begin{array}{cc}
0 & 1 \\
-q & 0
\end{array}
\right)
$$
$$
\epsilon _q{\cal D}_q=T^{+}\epsilon _qT=T\epsilon _qT^{+}
\eqno(2.8)
$$

The coproduct $\Delta $, counit $e$ and antipode $S$ on ${\cal D}_q$ defined
by

$$
\begin{array}{ccc}
\Delta ({\cal D}_q)={\cal D}_q\otimes {\cal D}_q, &
e({\cal D}_q)=1,&
S({\cal D}_q)={\cal D}_q^{-1}
\end{array}
\eqno(2.9)
$$

Adding ${\cal D}_q^{-1}$ to the algebra, we suppose that ${\cal D}_q\neq 0$.
Using (2.4,2.5) to obtain the commutation relations between group elements
$$
T=\left(
\begin{array}{cc}
a & b \\
c & d
\end{array}
\right) =\left(
\begin{array}{cc}
T_1^1 & T_2^1 \\
T_1^2 & T_2^2
\end{array}
\right)
\eqno(2.10)
$$
it can obtain the following component relations [3]
$$
\begin{array}{lll}
ab=q\;ba & bc=cb & cd=q\;dc \\
ac=q\;ca & bd=q\;db &  \\
ad=da+\lambda \;bc, &  &
\end{array}
\eqno(2.11)
$$

The quantum determinant ${\cal D}_q$ is central element and commutes with
all group elements
$$
T^i_k{\cal D}_q= {\cal D}_qT^i_k
\eqno(2.12)
$$

One can introduce the left-hand side exterior derivative
$$
\delta _Lf=\left( f{\frac{\stackrel{\leftarrow }{\partial }}{\partial T^i_k}}
\right) \delta _LT^i_k
\eqno(2.13)
$$
satisfying 1) Leibnitz rule

$$
\delta _L(f\theta )=(f\theta )\stackrel{\leftarrow }{\delta }_L=f(\delta
_L\theta )+(-1)^{\eta(\theta)}(\delta _Lf)\theta
\eqno(2.14)
$$
and 2) Poincare lemma
$$
\delta _L^2f=0
\eqno(2.15)
$$
where $\eta(\theta)$ is the degree of differentional form $\theta $ and
$f$ is arbitrary function of the group elements.

We wish to construct a quantum algebra of the left vector fields on the
quantum group $GL_q(2,C)$. To do this we need to introduce an infinitely
small neighborhood $\delta _LT$ of the unity of the group, and we need to
determine the commutation relations between the parameters of the group and
its differentials. Equivalently, we need to determine the commutation
relations between the parameters of the group and the Maurer-Cartan left
1-forms, which are given by
$$
\theta =S(T)\delta _LT=\left(
\begin{array}{cc}
\theta ^1 & \theta ^2 \\
\theta ^3 & \theta ^4
\end{array}
\right) =\left(
\begin{array}{cc}
\theta _1^1 & \theta _2^1 \\
\theta _1^2 & \theta _2^2
\end{array}
\right)
\eqno(2.16)
$$
where

$$
\begin{array}{ll}
\theta ^1= {\cal D}_q^{-1}(d\delta _La-q^{-1}b\delta _Lc)=\theta
_1^1~~~~~~ & \theta ^3=D_q^{-1}(a\delta _Lc-qc\delta _La)=\theta_1^2 \\
&\\
\theta ^2= {\cal D}_q^{-1}(d\delta _Lb-q^{-1}b\delta _Ld)=\theta _2^1
&
\theta ^4=D_q^{-1}(a\delta _Ld-qc\delta _Lb)=\theta_2^2
\end{array}
\eqno(2.17)
$$

Left coaction $\Delta _L$ extends to Maurer-Cartan left 1-forms such that
$$
\Delta _L\circ\delta _L = (id\otimes \delta _L)\circ\Delta
$$
$$
\Delta _L(\delta _LT)=\left(
\begin{array}{cc}
a & b \\
c & d
\end{array}
\right) \otimes \left(
\begin{array}{cc}
\delta a & \delta b \\
\delta c & \delta d
\end{array}
\right)
\eqno(2.18)
$$
\medskip
$$
\Delta _L(\theta (\delta _L))=1\otimes \theta (\delta _L)
$$

In the same way we define the right-hand side exterior derivative
$$
\stackrel{\rightarrow }{\delta _R}f={\stackrel{\rightarrow }{\delta
_RT_i^k}}\;\; {\frac{\stackrel{\rightarrow }{\partial }}{\partial T_i^k}}f
\eqno(2.19)
$$
satisfying the classical Leibnitz rule
$$
\stackrel{\rightarrow }{\delta _R}(\omega g)=(\stackrel{\rightarrow }{\delta
_R}\omega )g+(-1)^{\eta(\omega)}\omega (\stackrel{\rightarrow }{\delta
_R}g)
\eqno(2.20)
$$
and Poincare lemma $\stackrel{\rightarrow }{\delta _R}^2f=0$.

Let us define the Maurer-Cartan right-invariant 1-form $\omega $ on $T$
$$
\omega =\delta _RT\;S(T)=\left(
\begin{array}{cc}
\omega ^1 & \omega ^2 \\
\omega ^3 & \omega ^4
\end{array}
\right) =\left(
\begin{array}{cc}
\omega _1^1 & \omega _2^1 \\
\omega _1^2 & \omega _2^2
\end{array}
\right) ,
\eqno(2.21)
$$
where

$$
\begin{array}{ll}
\omega ^1=(
\stackrel{\rightarrow }{\delta _Ra}d-q\;\stackrel{\rightarrow }{\delta
_Rb}c){\cal D}_q^{-1}=\omega _1^2,~~~~~~
&
\omega ^3=(
\stackrel{\rightarrow }{\delta _Rc}d-q\;\stackrel{\rightarrow }{\delta
_Rd}c){\cal D}_q^{-1}=\omega _1^2, \\
&\\
\omega^2 =(\stackrel{\rightarrow }{ \delta _Rb}a-\frac
1q\;\stackrel{\rightarrow }{\delta _Ra}b\;){\cal D} _q^{-1}=\omega _2^1,
&
\omega ^4=(\stackrel{\rightarrow }{
\delta _Rd}a-\frac 1q\stackrel{\rightarrow }{\delta _Rc}b){\cal D}
_q^{-1}=\omega _2^2,
\end{array}
\eqno(2.22)
$$

We note, that left differential arranges after left derivatives, but right
differential before right derivatives.

Right coaction $\Delta_R (\stackrel{\rightarrow}{\delta_R}) = \stackrel{
\rightarrow}{\delta_R}\otimes 1$:
$$
\Delta_R(\stackrel{\rightarrow}{\delta_R}T)=\left(
\begin{array}{cc}
\stackrel{\rightarrow}{\delta_R}a & \stackrel{\rightarrow}{\delta_R}b \\
\stackrel{\rightarrow}{\delta_R}c & \stackrel{\rightarrow}{\delta_R}d
\end{array}
\right) \otimes \left(
\begin{array}{cc}
a & b \\
c & d
\end{array}
\right)
\eqno(2.23)
$$
$$
\Delta_R(\omega(\stackrel{\rightarrow}{\delta_R})) = \omega\otimes 1
$$

Notice, that on the group parameters $T^i_k$ the action $\Delta_R$ and $
\Delta_L$ coincides with the action $\Delta$.

To specify a differential calculus it is necessary to define commutation
relations between $T^i_k$ and their differentials $\delta T^i_k$.

Throughout the recent development of differential calculus on quantum groups
and quantum spaces, two principal concepts are readily seen. First of them,
for\-mu\-la\-ted by Wo\-ro\-no\-wicz [5], is known as bicovariant
differential calculus on quantum groups:
$$
(1\otimes\Delta_R)\Delta_L = (\Delta_L\otimes 1)\Delta_R
\eqno(2.24)
$$

The examples were considered by Schupp, Wats, Zumino [11],
Jur\v co[8], Sudberry [12], M\"uller-Hoissen [13,14], Isaev,~Pjatov [15],
Demidov [16], Manin [7] and others [17--20]. But the natural way
of obtaining the $SL_q(N)$-differential calculus by performing reduction
from the $GL_q(N)$-calculus (${\cal D}_q= 1$) failed in the quantum case.

Another concept, introduced by Woronowicz [4] and
Schirrmacher, Wess, Zumino [10] proceeds from the requirement
of left- (right-) invariant differential calculus only. We will consider
the last concept.

In this paper we propose the construction of the differential calculus on the
group and its subgroup with the property of ''quantum matrjoshka'' (matched
reduction): the differential calculus on the quantum group $GL_q(2,C)$ has to
contain the differential calculus on the quantum subgroup $SL_q(2,C)$ and
quantum plane $C_q(2|0)$ .

We postulate, that the quantum determinant ${\cal D}_q$ is a central element
for $\theta $ or $\omega $ 1-forms
$$
\theta ^i{\cal D}_q={\cal D}_q\theta ^i,D_q\omega =\omega D_q
\eqno(2.25)
$$
in contrast to [11,16,18], where $\theta ^i{\cal
D}_q=q^{-2} {\cal D}_q\theta ^i$ for the bicovariant differential calculus.
A consequence of these conditions is that we can take ${\cal
D}_q=\det _qT=1$ for $SL_q(2,C)$ and will obtain the matched differential
calcules on the $ GL_q(2,C)$ and the $SL_q(2,C)$.

\section{Left differentioal calculus on the $GL_q(2,C)$.}

Let us first consider the case of the left invariant differential calculi
on the $GL_q(2,C)$. Matched differential calculus has to satisfy seven
basic conditions, six of them are the conditions of consistency with RTT
relations and last condition is the condition of consistency of
differential calculus on the quantum group and its subgroup.

\setcounter{nnn}{1}\Roman{nnn}) Commutation relations between the $T^i_k$
and their differentials can be expressed in terms of the Maurer--Cartan
left 1-forms[13]

$$
\theta^{\alpha}_{\beta}T^i_k = T^i_m \left( \bara^m_k\right)
^{\alpha\gamma}_{\delta\beta} \theta^{\delta}_{\gamma}
\eqno(3.1)
$$
where $\bara^m_k = \lmatr{cc} A & B\\ C & D \rmatr$ and $A,B,C,D$ are
$4\times 4$ matrices with complex entries.
\medskip
\medskip

\setcounter{nnn}{2}\Roman{nnn}) The requirement, that
$\theta^{\alpha}_{\beta}$ commute with commutation relations (2.4)
$$
\theta^i_j\left(R^{\alpha n}_{\beta m} T^{\beta}_{\delta} T^m_k -
T^n_m T^{\alpha}_{\beta} R^{\beta m}_{\delta k}\right) = 0
\eqno(3.2)
$$
leads to some additional equations for the matrices $A,B,C,D$
$$
R^{\alpha n}_{\beta m} \bara^{\beta}_{\delta} \bara^m_k =
\bara^n_m \bara^{\alpha}_{\beta} R^{\beta m}_{\delta k}
\eqno(3.3)
$$
$A,B,C,D$ have to form a representation of $a,b,c,d$
$$
\begin{array}{c}
\begin{array}{lll}
AB=q\;BA & AC=q\;CA & BC=CB \\
BD=q\;DB & CD=q\;DC & \\
\end{array}\\
AD=DA+(q-\frac{1}{q})BC;
\end{array}
\eqno(3.4)
$$

{}From the definition of the left $\theta$-form, as $\theta= S(T^{-1})
\delta_L T$ we find
$$
\delta_L T^{\beta}_{\alpha}=T^{\gamma}_{\alpha}\theta^{\beta}_{\gamma}
\eqno(3.5)
$$
and now we can use left exterior derivative $\delta_L$.

\setcounter{nnn}{3}\Roman{nnn}) For the
quantum determinant $\qdet$ we define quantum trace by the definition

$$
\delta_L\qdet=\qdet\qtr\theta
\eqno(3.6)
$$
Using (3.1) and (3.5) we obtain
$$
\begin{array}{rcl}
\delta_L(ad\;-\;q\;bc) &=& \qdet\left(\theta^1+A^4_l\theta^l-
	\frac{1}{q}B^3_l\theta^l\right) =
	\qdet\left(\theta^4+D^1_l\theta^l - qC^2_l\theta^l\right)\\
&&\\
\delta_L(ad\;-\;q\;cb) &=& \qdet\left(\theta^4+D^1_l\theta^l-
	qB^3_l\theta^l\right)\\
&&\\
\delta_L(da\;-\;\frac{1}{q}\;bc) &=&
\qdet\left(\theta^1+A^4_l\theta^l- \frac{1}{q}C^2_l\theta^l\right)\\
&&\\
\delta_L(da\;-\;\frac{1}{q}\;cb) &=&
\qdet\left(\theta^1+A^4_l\theta^l- \frac{1}{q}B^3_l\theta^l\right) =
	\qdet\left(\theta^4+D^1_l\theta^l - qC^2_l\theta^l\right)
\end{array}
\eqno(3.7)
$$
\medskip
and additional equations
$$
\begin{array}{rcl}
B^3_l &=& C^2_l\\
&&\\
\theta^2+B^1_l\theta^l-qA^2_l\theta^l &=& 0\\
&&\\
\theta^3+C^4_l\theta^l-\frac 1q D^3_l\theta^l &=& 0
\end{array}
\eqno(3.8)
$$

\setcounter{nnn}{4}\Roman{nnn}) Using the Maurer--Cartan equations
$\delta\theta = -\theta\land\theta$ we find the conditions of consistency
with RTT relations
$$
\delta_L(R_{12}T_1T_2 - T_2T_1R_{12}) = 0
\eqno(3.9)
$$
as a consequence we obtain the same equations as (2.27--2.28).

\setcounter{nnn}{5}\Roman{nnn}) Using (3.2) we receive
$$
\theta^i_k\qdet = \qdet\left(AD-qBC\right)^i_m\theta^m_k
$$
and as a consequence we obtain (3.4).

\setcounter{nnn}{6}\Roman{nnn}) Differentiating commutation relations
(3.1) and using the Maurer--Cartan equation
we find some additional equations for $A,B,C,D$
$$
\delta_L\theta^{\alpha}_{\beta}T^n_k - \theta^{\alpha}_{\beta} \delta_L
T^n_k
=
\delta_LT^n_m\left(\bara^m_k\right)^{\alpha\gamma}_{\delta\beta}
\theta^{\delta}_{\gamma} + T^n_m
\left(\bara^m_k\right)^{\alpha\gamma}_{\delta\beta} \delta_L
\theta^{\delta}_{\gamma} - \theta^{\alpha}_{\gamma}
\theta^{\gamma}_{\beta} T^n_k - \theta^{\alpha}_{\beta} T^n_{\gamma}
\theta^{\gamma}_k =
$$
$$
=
T^n_{\rho}\theta^{\rho}_m \left(\bara^m_k\right)^{\alpha\gamma}_{
\delta\beta}\theta^{\delta}_{\gamma} -
T^n_m\left(\bara^m_k\right)^{\alpha\gamma}_{\delta\beta}
\theta^{\delta}_{\rho} \theta^{rho}_{\gamma}
$$
or
$$
T^n_f\left[\left(\bara^f_m\right)^{\alpha\varphi}_{\tau\gamma}
\left(\bara^m_k\right)^{\gamma\beta}_{\sigma\beta} +
\left(\bara^f_{\sigma}\right)^{\alpha\varphi}_{\tau\beta}
\delta^{\rho}_k +
\left(\bara^{\varphi}_k\right)^{\alpha\rho}_{\gamma\beta}
\delta^f_{\tau} -
\left(\bara^f_k\right)^{\alpha\rho}_{\tau\beta}
\delta^{\sigma}_{\varphi}\right]
\theta^{\tau}_{\varphi}\theta^{\sigma}_{\rho} = 0
$$

\medskip
We have $q$-deformation of algebra of $\theta$-forms:
$$
\left(\bara^f_m\right)^{\alpha\varphi}_{\tau\gamma}
\left(\bara^m_k\right)^{\gamma\rho}_{\sigma\beta}
\theta^{\tau}_{\varphi} \theta^{\sigma}_{\rho} +
\left(\bara^f_{\sigma}\right)^{\alpha\varphi}_{\tau\beta}
\theta^{\tau}_{\varphi} \theta^{\sigma}_{k} +
\left(\bara^{\varphi}_{k}\right)^{\alpha\rho}_{\sigma\beta}
\theta^{f}_{\varphi} \theta^{\sigma}_{\rho} -
\left(\bara^{f}_{k}\right)^{\alpha\rho}_{\tau\beta}
\theta^{\tau}_{\sigma} \theta^{\sigma}_{\rho}
= 0
\eqno(3.10)
$$

\setcounter{nnn}{7}\Roman{nnn}) The last condition is the basic condition
which differs matched differential calculus from the bicovariant
differential calculus
$$
\theta^k\qdet = \qdet\theta^k
$$
which leads to the equation $AD - qBC = 1$.

To solve the system of equations for the matrix $\bara$ following the
conditions
(\setcounter{nnn}{1}\Roman{nnn}--\setcounter{nnn}{7}\Roman{nnn}) we will
consider the representation for entries of $\bara$ with
$$
B=C=0\;\;\;\;\;\;AD=1
$$
{}From (3.7) we have obtained only two different expressions
$$
\delta_L\qdet = \qdet\left[D^1_1 \theta^1 +
(1+D^1_4)\theta^4\right],
\eqno(3.11)
$$
$$
\delta_L\qdet =
\qdet\left[(1+A^4_1)\theta^1 + A^4_4 \theta^4\right],
$$
where we used composite indeces 1,2,3,4.

Finally we find matrices $A,D$ in terms of independent variables
$\beta=D^1_1,\;\alpha=A^4_4$
\medskip
$$
\begin{array}{c}
A=\lmatr{cccc} 1-\alpha+\frac{\alpha}{\beta} & 0 & 0 &
\frac{(1-\alpha)\alpha}{\beta}\\
0	 & \frac{1}{q}& 0 & 0 \\
0	 & 0     & \frac{1}{q} & 0\\
\beta-1 & 0 & 0 & \alpha\\
\rmatr \;\;,\;\;\;\;
D=\lmatr{cccc}
\beta & 0 & 0 & \alpha-1\\
0 & q & 0 & 0\\
0 & 0 & q & 0\\
\frac{(1-\beta)\beta}{\alpha} & 0 & 0 &
1-\beta+\frac{\beta}{\alpha}
\rmatr
\end{array}
\eqno(3.12)
$$
\medskip

{}From the commutation relations for $\theta$-forms for given $A$ and
$D$ we obtain the following commutation relations
$$
\theta^1\theta^4 + \theta^4\theta^1 = -\frac{(\alpha-1)}{\beta}
\left(\theta^4\right)^2 - \frac{(\beta-1)}{\alpha} \left(\theta^1\right)^2
- \frac{-\beta+q^2\alpha}{\beta\alpha}\theta^2\theta^3 =
\eqno(3.13)
$$
\nopagebreak
$$
= -\frac{(\beta+\alpha-
\beta\alpha)}{(1-\beta)\beta} \left(\theta^4\right)^2 -
\frac{(\beta-\alpha-\beta^2)}{\beta+2\alpha- \beta\alpha}
\left(\theta^1\right)^2 -
\frac{(-\beta+q^2\alpha)(1-q^2-\alpha)\alpha}
{(\beta+2\alpha-\beta\alpha)(1-\beta)\beta}\theta^2\theta^3 =
$$
\nopagebreak
$$
= -\frac{(-\beta+\alpha-\alpha^2)}
{\alpha+2\beta-\beta\alpha}\left(\theta^4\right)^2 -
\frac{(\beta+\alpha-\beta\alpha)}
{(1-\alpha)\alpha}\left(\theta^1\right)^2  -
\frac{(-\beta+q^2\alpha)(1-q^{-2}-\alpha)\beta}
{(1-\beta)(\alpha+2\beta-\beta\alpha)\alpha} \theta^2\theta^3
$$
\nopagebreak
$$
(\theta^2)^2=(\theta^3)^2=0,\;\;\;\;\;\theta^2\theta^3+q^2\theta^3\theta^2=0
\eqno(3.14)
$$
which are compatible, when
$$
\left(\theta^1\right)^2=0;\;\;
\left(\theta^4\right)^2=0;\;\;
\beta=q^2\alpha
\eqno(3.15)
$$
Using (3.15), we have
$$
\theta^1\theta^4 + \theta^4\theta^1 =0
\eqno(3.16)
$$
The remaining commutation relations are
$$
\begin{array}{lcl}
q^2\alpha\theta^2\theta^1 + (q^{-2}+1-\alpha)\theta^1\theta^2 +
(\alpha-1)\theta^2\theta^4 + q^{-2}(1-\alpha)\theta^4\theta^2 &=&0\\
&&\\
q^2(1-q^2\alpha)\theta^2\theta^1 + (q^2\alpha-1)\theta^1\theta^2 +
(1+q^2-q^2\alpha)\theta^2\theta^4 + \alpha\theta^4\theta^2 &=&0\\
&&\\
q^2\alpha\theta^1\theta^3 + (q^{-2}+1-\alpha)\theta^3\theta^1 +
(\alpha-1)\theta^4\theta^3 + q^{-2}(1-\alpha)\theta^3\theta^4 &=&0\\
&&\\
q^2(1-q^2\alpha)\theta^1\theta^3 + (q^2\alpha-1)\theta^3\theta^1 +
(1+q^2-q^2\alpha)\theta^4\theta^3 + \alpha\theta^3\theta^4 &=&0
\end{array}
\eqno(3.17)
$$

In fact the parameter $\alpha$ plays the role of the definition of the
$q$-trace [3,18,21]
$$
\delta_L\qdet=\alpha\qdet(q^2\theta^1+\theta^4)
$$
Thus we obtained the one--parameter family differential calculus on the
$GL_q(2,C)$.

The next step is to construct a left vec\-tor field al\-geb\-ra for the
quan\-tum
group \linebreak $GL_q(2,C)$. By definition, the effect of applying the left
differential to an arbitrary function on the quantum group is

$$
\delta_Lf= \frac{\lpar f}{\partial T^i_k}\delta_L T^i_k = \left(f
\anabla_k\right)\theta^k
\eqno(3.18)
$$
where $\anabla_k$ are the left vector fields on the quantum group.

{}From the lemma Poincare
$$
\delta_L^2f=-(f\anabla_k)\delta_L\theta^k-(f\anabla_k)(\anabla_l\theta^l)
\theta^k=0
$$
and the Maurer--Cartan equations we find the algebra of the left
vector fields for an arbitrary $\alpha$
$$
\anabla_3\anabla_2 -
q^2\anabla_2\anabla_3 = \hanabla_1,
$$
$$
q^2\anabla_2\hanabla_1 - q^{-2}\hanabla_1\anabla_2 = (1+q^2)\anabla_2
$$
$$
q^2\anabla_1\anabla_3 - q^{-2}\anabla_3\anabla_1 = (1+q^2)\anabla_3
\eqno(3.19)
$$
$$
\left[\anabla_4,\anabla_2\right] = -(1+q^2)(1-\alpha)\anabla_2\hanabla_1
+ [q+(1+q^2)(1-\alpha)]\anabla_2
$$
$$
\left[\anabla_3,\anabla_4\right] = -(1+q^2)(1-\alpha) \hanabla_1\anabla_3
+ [1+(1+q^2)(1-\alpha)]\anabla_3
$$
$$
\left[\anabla_4,\hanabla_1\right] = 0
$$
where $\hanabla_1$ represents the combination
$$
\hanabla_1=\anabla_1-q^2\anabla_4
\eqno(3.20)
$$

In paper [22] authors have investigated the case of differential
calculus with $\alpha=1$.

Now we find, that the algebra of 1--forms (3.17) and the algebra of vector
fields can be decomposed on the algebra $SL_q(2,C)$ and
$U(1)$ subalgebra only for unique value of the parameters
$$
\alpha=\frac{2}{1+q^2}, \;\;\;\;\;\; \beta=\frac{2q^2}{1+q^2}.
\eqno(3.21)
$$

In this case the commutation relations between the parameters and 1--forms
are diagonalized after the choice of the new basis of 1--forms
$$
\ttheta^1=\frac{2}{q+1/q}(q\theta^1+\frac{1}{q}\theta^4)=\qtr\theta
$$
$$
\ttheta^4=\frac{1}{1+q^2}(\theta^1-\theta^4)
$$
$$
[\ttheta^1,a]=[\ttheta^1,d]=0
\eqno(3.22)
$$
$$
\begin{array}{lcl}
\ttheta^4a=q^{-2}a\ttheta^4 &~~~~~& \ttheta^4d=q^2d\ttheta^4,\\
\theta^2a=q^{-1}a\theta^2 && \theta^2d=qd\theta^2,\\
\theta^3a=q^{-1}a\theta^3 && \theta^3d=qd\theta^3
\end{array}
$$

The other relations are found by making the interchanges $a\leftrightarrow
c,d\leftrightarrow b$.

The commutation relations (3.17) can be written as
$$
\begin{array}{ll}
\ttheta^1\theta^2 + \theta^2\ttheta^1 =0,~~~~~~ &
\ttheta^1\theta^3 + \theta^3\ttheta^1 =0\\
q^2\theta^2\ttheta^4 + q^{-2}\ttheta^4\theta^2 =0, &
q^2\ttheta^4\theta^3 + q^{-2}\theta^3\ttheta^4 =0\\
\ttheta^1\ttheta^4 + \ttheta^4\ttheta^1 =0, &
\theta^2\theta^3 + q^{-2}\theta^3\theta^2 =0
\end{array}
\eqno(3.23)
$$
$$
(\ttheta^1)^2=(\theta^2)^2=(\theta^3)^2=(\ttheta^4)^2=0
$$

For the consistency of the differential algebra of the 1-forms we should
define how to order lexicographycally any cubic monomial (Diamond
Lemma[7]). For example, if we try to express the
cubic monomial $\ttheta^4\theta^3\theta^2$ of the permutations ($432\to
423\to 243\to 234$) and of the permutations ($432\to 342\to 324\to 234$),
so we should receive the same result
$$
\ttheta^4\theta^3\theta^2\to
-q^{-4}\theta^3\ttheta^4\theta^2\to \theta^3\theta^2\ttheta^4\to
-q^2\theta^2\theta^3\ttheta^4;
$$
$$
\ttheta^4\theta^3\theta^2\to
-q^2\ttheta^4\theta^2\theta^3\to q^6\theta^2\ttheta^4\theta^3\to
-q^2\theta^2\theta^3\ttheta^4
$$

An attempt to order the monomials $\ttheta^4\theta^2\ttheta^1$ and
$\ttheta^4\theta^3\ttheta^1$ using (3.19--3.20) leads to
$$
\left\{
\begin{array}{l}
\ttheta^4\theta^2\ttheta^1\to -\ttheta^4\ttheta^1\theta^2\to
\ttheta^1\ttheta^4\theta^2\to -q^4\ttheta^1\theta^2\ttheta^4\\
\ttheta^4\theta^2\ttheta^1\to -q^4\theta^2\ttheta^4\ttheta^1\to
q^4\theta^2\ttheta^1\ttheta^4\to -q^4\ttheta^1\theta^2\ttheta^4
\end{array}
\right.
$$
and
$$
\left\{
\begin{array}{l}
\ttheta^4\theta^3\ttheta^1\to -\ttheta^4\ttheta^1\theta^3\to
\ttheta^1\ttheta^4\theta^3\to -q^{-4}\ttheta^1\theta^3\ttheta^4\\
\ttheta^4\theta^3\ttheta^1\to -q^{-4}\theta^3\ttheta^4\ttheta^1\to
q^{-4}\theta^3\ttheta^1\ttheta^4\to -q^{-4}\ttheta^1\theta^3\ttheta^4
\end{array}
\right.
$$

Normal ordering of monomials like $\ttheta^4\theta^3\ttheta^4$,
$\theta^2\ttheta^4\theta^2$, $\theta^3\ttheta^1\theta^3$,
$\ttheta^1\theta^2\ttheta^1$ and so on are easily shown to vanish.

Thus, the different ways of ordering cubic monomial lead to the same result
and it guarantees the consistency of the differential algebra $GL_q(2,C)$.

In conclusion of this section we find the commutation relations between
the parameters of the quantum group and its differentials.

For any parameters of the group $K$ and $N$ we may decompose the
expression $K\delta N$ in the complete basis of 1--forms
$$
\begin{array}{rcl}
K\delta N &=& \tilde A \delta a + \tilde B \delta b + \tilde C\delta c +
\tilde D \delta d =\\
&=& (\tilde Aa+\tilde Cc)\theta^1 + (\tilde Ba+\tilde Dc)\theta^2 +
(\tilde Ab+\tilde Cd)\theta^3 + (\tilde Bb+\tilde Dd)\theta^4
\end{array}
$$
where $\tilde A$,$\tilde B$,$\tilde C$,$\tilde D$ are arbitrary functions
of parameters of quantum group.

Using the commutation relations between the parameters and 1--forms we
obtain the commutation relations between the parameters and its
differentials.
$$
\begin{array}{l}
\delta_L a\;a=q^{-2}\;a\;\delta_L a +
\frac{q^2-1}{2q^2}\;a^2\;\mbox{Tr}_q\theta,
\\
\delta_L c\;c=q^{-2}\;c\;\delta_L c +
\frac{q^2-1}{2q^2}\;c^2\;\mbox{Tr}_q\theta,
\\
\delta_L a\;c=q^{-1}\;c\;\delta_L a +
\frac{q^2-1}{2q^2}\;a\;c\;\mbox{Tr}_q\theta,
\\
\delta_L c\;a=q^{-1}\;a\;\delta_L c + (q^{-2}-1)\;c\;\delta_L a +
\frac{q^2-1}{2q^2}\;c\;a\;\mbox{Tr}_q\theta,
\\
\delta_L b\;b=q^{2}\;b\;\delta_L b +
\frac{1-q^2}{2}\;b^2\;\qtr\theta,
\\
\delta_L d\;d=q^{2}\;d\;\delta_L d +
\frac{1-q^2}{2}\;d^2\;\qtr\theta,
\\
\delta_L b\;d=q\;d\;\delta_L b + (q^2-1)\;b\;\delta_L d +
\frac{1-q^2}{2}\;b\;d\;\qtr\theta,
\\
\delta_L d\;b=q\;b\;\delta_L d +
\frac{1-q^2}{2}\;d\;b\;\qtr\theta,
\end{array}
$$
$$
\begin{array}{l}
\delta_L a\;b=q\;b\;\delta_L a +
\frac{(q^2-1)}{q^2\qdet}\;a\;b\;(q\;c\;\delta_L b - a\;\delta_L d) +
\frac{q^2-1}{2q^2}\;a\;b\;\qtr\theta,
\\
\delta_L a\;d=d\;\delta_L a +
\lambda\;b\;\delta_L c +
\frac{q^2-1}{\qdet}\;a\;d(d\;\delta_L a - \frac 1q\;b\;\delta_L c) -
\frac{q^2-1}{2}\;a\;d\;\qtr\theta
\\
\delta_L c\;b=b\;\delta_L c + \frac{(q^2-1)}{\qdet}\;c\;b\; (d\;\delta_L a -
\frac 1q\;b\;\delta_L c) - \frac{(q^2-1)}{2}\;c\;b\;\qtr\theta,
\\
\delta_L c\;d=q\;d\;\delta_L c + \frac{(q^2-1)}{\qdet}\;c\;d\; (d\;\delta_L a -
\frac 1q\;b\;\delta_L c) - \frac{(q^2-1)}{2}\;c\;d\;\qtr\theta,
\\
\delta_L b\;a=q^{-1}\;a\;\delta_L b +
\frac{(q^2-1)}{q^2\qdet}\;b\;a\;(q\;c\;\delta_L b - a\;\delta_L d) +
\frac{(q^2-1)}{2q^2}\;b\;a\;\qtr\theta,
\\
\delta_L b\;c=c\;\delta_L b + \frac{(q^2-1)}{\qdet}\;b\;c\; (d\;\delta_L a -
\frac 1q\;b\;\delta_L c) - \frac{(q^2-1)}{2}\;b\;c\;\qtr\theta,
\\
\delta_L d\;a=a\;\delta_L d - \lambda\;c\;\delta_L b + \frac{q^2-1}{\qdet}\;
d\;a\;(d\;\delta_L a - \frac 1q\;b\;\delta_L c) -
\frac{(q^2-1)}{2}\;d\;a\;\qtr\theta
\\
\delta_L d\;c=q^{-1}\;c\;\delta_L d + \frac{(q^2-1)}{\qdet}\;d\;c\;
(d\;\delta_L a - \frac 1q\;b\;\delta_L c) -
\frac{(q^2-1)}{2}\;d\;c\;\qtr\theta,
\end{array}
\eqno(3.24)
$$
where $\qtr\theta=\frac{2}{q+1/q}(q\theta^1+\frac 1q \theta^4)=\ttheta^1$.

The algebra of vector fields in this case has the following form
$$
\anabla_3\anabla_2-q^2\anabla_2\anabla_3=\hanabla_1
$$
$$
q^2\anabla_2\hanabla_1-q^{-2}\hanabla_1\anabla_2=(1+q^2)\anabla_2
\eqno(3.25)
$$
$$
q^2\hanabla_1\anabla_3-q^{-2}\anabla_3\hanabla_1=(1+q^2)\anabla_3
$$
$$
[\hanabla_4,\hanabla_1]=[\hanabla_4,\anabla_2]=[\hanabla_4,\anabla_3]=0
$$
where
$$
\hanabla_1 = \anabla_1-q^2\anabla_4,~~~~~~~
\hanabla_4 = \anabla_1 + \anabla_4
\eqno(3.26)
$$

After the mapping $\hanabla_1,\anabla_2,\anabla_3,\hanabla_4$ to the new
generators $H,T_2,T_3,N$
$$
\hanabla_1=\frac{1-q^{-2H}}{1-q^{-2}}\;\;\;,\;\;\;
\anabla_2=q^{-H/2}T_2,
$$
\nopagebreak
$$
\hanabla_4=N\;\;\;\;\;\;,\;\;\;\;\;\;
\anabla_3=q^{-H/2}T_3,
\eqno(3.27)
$$
we obtain the algebra $U_qGL(2,C)$ in the form of a Drinfeld--Jimbo
algebra [1,2]
$$
[T_3,T_2]=\frac{q^H-q^{-H}}{q-q^{-1}}
$$
$$
[H,T_3]=2T_3\;\;\;\;,\;\;\;
[H,T_2]=-2T_2
\eqno(3.28)
$$
$$
[N,H]=[N,T_2]=[N,T_3]=0
$$

Note, that there exists one more matched solution when $\alpha\to 0$ and
$\beta\to 0$, but its ratio remains equal $\frac{\beta}{\alpha}=q^2$.
In this case we have from (3.12)

$$
\begin{array}{cc}
\begin{array}{rcl}
(\theta^1+\theta^4)a &=& \frac{1}{q^2}a(\theta^1+\theta^4)\\
(\theta^1+q^2\theta^4)a &=& a(\theta^1+q^2\theta^4)\\
\theta^2a &=& \frac 1q a\theta^2\\
\theta^3a &=& \frac 1q a\theta^3
\end{array} &
\begin{array}{rcl}
\theta^1d&=&-d\;\theta^4\\
\theta^4\;d&=&(1+q^2)d\theta^4+q^2d\theta^1\\
\theta^1a&=&(1+q^{-2})a\theta^1+q^{-2}a\theta^4\\
\theta^4a&=&-a\theta^1
\end{array}\\
&\\
\begin{array}{rcl}
(\theta^1+\theta^4)d &=& q^2d(\theta^1+\theta^4)\\
(\theta^1+q^2\theta^4)d &=& d(\theta^1+q^2\theta^4)\\
\theta^2d &=& q d\theta^2\\
\theta^3d &=& q d\theta^3
\end{array} &
\end{array}
\eqno(3.36)
$$

The other relations are found by making interchanges $d\to b$, $a\to c$.
The commutation relations between $\theta^2$,$\theta^3$ and
$a$,$b$,$c$,$d$ did not change.

\section{Left differential calculus on $SL_q(2,C)$ and quantum plane.}

To obtain the differential calculus on the $SL_q(2,C)$ it is necessary to
suppose additional constraints $\qdet=1$ and
$\delta\qdet=\qdet\qtr=0$. In an opposite way to the bicovariant calculus
it is possible, so the quantum determinant $\qdet$ commutes with 1--forms
and $\delta\qdet=0$ satisfies by vanishing $\ttheta^1=\qtr\theta =
\alpha(q\theta^1+\frac 1q \theta^4)=0$. As a consequence of these
conditions we obtained three--dimensional differential calculus, which is
independent of the parameters $\alpha$ and $\beta$.

Thus we have next commutation relations in the $SL_q(2,C)$ case
$$
\begin{array}{l}
\theta^1a=q^{-2}a\theta^1\\
\theta^2a=q^{-1}a\theta^2\\
\theta^3a=q^{-1}a\theta^3
\end{array}
\;\;\;,\;\;\;
\begin{array}{l}
\theta^1d=q^2d\theta^1\\
\theta^2d=qd\theta^2\\
\theta^3d=qd\theta^3
\end{array}
\eqno(4.1)
$$
$$
\begin{array}{l}
(\theta^1)^2=(\theta^2)^2=(\theta^3)^2=0,\;\;\;\;\;\;
\theta^4=-q^2\theta^1\\
\theta^1\theta^2+q^4\theta^2\theta^1=0\\
\theta^1\theta^3+q^{-4}\theta^3\theta^1=0\\
\theta^2\theta^3+q^{-2}\theta^3\theta^2=0
\end{array}
\eqno(4.2)
$$
and other relations are made by interchanges $a\leftrightarrow c$,
$d\leftrightarrow b$.

Notice, that similary commutation relations for $SU_q(2)$ were received by
\linebreak
Woronowich~[1]. The algebra of vector fields in this case has the
following form
$$
\begin{array}{rcl}
q^2\nabla_1\nabla_3-q^{-2}\nabla_3\nabla_1 &=& (1+q^2)\nabla_3\\
q^2\nabla_2\nabla_1-q^{-2}\nabla_1\nabla_2 &=& (1+q^2)\nabla_2\\
\nabla_3\nabla_2-q^2\nabla_2\nabla_3 &=& \nabla_1
\end{array}
\eqno(4.3)
$$

The commutation relations between the parameters and its differentials
have the form:
$$
\begin{array}{l}
\begin{array}{rcl}
\delta_L a\,a & = & q^{-2}a\,\delta_L a \\
\delta_L c\,c & = & q^{-2}c\,\delta_L c \\
\delta_L a\,c & = & q^{-1}c\,\delta_L a \\
\delta_L c\,a & = & q^{-1}a\,\delta_L c + (q^{-2}-1) c\,\delta_L a
\end{array}
\;\;\;,\;\;\;\;
\begin{array}{rcl}
\delta_L b\,b & = & q^2b\,\delta_L b \\
\delta_L d\,d & = & q^2d\,\delta_L d \\
\delta_L b\,d & = & q\,d\,\delta_L b +(q^2-1)b\,\delta_L d \\
\delta_L d\,b & = & q\,b\,\delta_L d
\end{array}
\\
\\
\begin{array}{rcl}
\delta_L a\,b & = & q\,b\,\delta_L a+(q^2-1)a\,b\,d\,\delta_L a + (
\frac{1}{q}-q) a\,b^2\delta_L c \\ \delta_L a\,d & = & q^2d\,\delta_L a +
q(q^2-1)b\,c\,d\,\delta_L a + (1-q^2)b^2c\,\delta_L c \\
\delta_L c\,b & = & b\,\delta_L c + (q^2-1)b\,c\,d\,\delta_L a +
(q^{-1}-q)b^2c\,\delta_L c \\
\delta_L c\,d & = & q\,d\,\delta_L c + (q^2-1)c\,d^2\delta_L a + (q^{-1}-q)
        c\,d\,b\,\delta_Lc \\
\delta_L b\,a & = & q^{-1}a\,\delta_L b + (q^2-1)b\,a\,d\,\delta_L a +
(1-q^2)b^2\,a\,\delta_L c \\
\delta_L b\,c & = & c\,\delta_L b + (q^2-1)b\,c\,d\,\delta_L a +
(q^{-1}-q)b^2c\,\delta_L c \\
\delta_L d\,a & = & q^{-2}a\,\delta_L d + q^{-2}(q^{-1}-q)b\,c\,a\,\delta_L d +
(1-q^{-2}) b\,c^2\delta_L b \\
\delta_L d\,c & = & q^{-1}c\,\delta_L d + (q^{-2}-1)d\,c\,a\,\delta_L d +
(q-q^{-1})d\,c^2\delta_L b
\end{array}
\end{array}
\eqno(4.4)
$$

Let us mention that for the differential calculus described in (4.4)
there are the quartic powers of the group elements essentially in contrast
to the commutation relations between $T^i_k$ and the Maurer--Cartan left
invariant 1--forms. The quartic powers in the commutation relations have
been obtained in papers [13,14].

Another difference, contrary to the bicovariant R-matrix approach of
Manin[7], Demidov[16], Castellani, Ashieri[17] appears. Let us call to
mind the differential calculus on the quantum plane $C_q(2|0)$ of
Wess--Zumino[6]:

$$
\begin{array}{llcll}
\mbox{\setcounter{nnn}{1}\Roman{nnn})} & xy = q\,yx &~~~~~~~&
\mbox{\setcounter{nnn}{2}\Roman{nnn})} & q\rightarrow
q^{-1}\;\;\;x\rightarrow y\;\;\; y\rightarrow x \\
& \delta x\,x = q^{-2}x\,\delta x &&
& \delta x\,y = q^{-1}\,y\,\delta x \\
& \delta y\,y = q^2y\,\delta y &&
& \delta y\,x = q^{-1}\,x\,\delta y + (q^{-2}-1)\,y\,\delta x\\
& \delta x\,y = q\,y\,\delta x + (q^2-1)x\,\delta y &&&\\
& \delta y\,x = q\,x\,\delta y &&&
\end{array}
\eqno(4.5)
$$

It is possible to introduce a quantum plane $C_q(2|0)$ as the Hopf algebra
surjection $\pi\,:SL_q(2,C)\to C_q(2|0)$ such that $\pi(T^1_2)=\pi(b)=0$ or
$\pi(T^2_1)=\pi(c)=0$:
$$
{\cal D}_q={\det}_qT=ad=1
$$
$$
T_-=\left(
\begin{array}{cc}
a & 0 \\
c & d
\end{array}
\right) = \left(
\begin{array}{cc}
x & 0 \\
y & x^{-1}
\end{array}
\right);
$$
\medskip
$$
T_+=\left(
\begin{array}{cc}
a & b \\
0 & d
\end{array}
\right) = \left(
\begin{array}{cc}
y^{-1} & x \\
0 & y
\end{array}
\right);
$$
\medskip

We see, that in the case $\pi(c)=0$ or $\pi(b)=0$ commutation relations (4.4)
form the commutation relations on the quantum plane.

The algebra of left 1--forms and the algebra of vector fields can be
obtained if we put in (4.2)
$$
\begin{array}{llll}
1) & c=0, & \theta^3=0, & \nabla_3=0\\
2) & b=0, & \theta^2=0, & \nabla_2=0
\end{array}
$$

Hence, we received the differential calculus on the quantum group
$GL_q(2,C)$ and its subgroup with the property of natural reduction
("quantum matrjoshka"): the differential calculus on the $GL_q(2,C)$
contains the differential calculus on the $SL_q(2,C)$, which also contains
the differential calculus on the quantum plane $C_q(2|0)$.

\section{Right differential calculus on $GL_q(2,C)$.}

In this section we give a construction of the right differential calculi
on $GL_q(2,C)$.

Analogously to previous sections we will assume, that the commutation
relations between the right 1-forms $\omega$ and parameters
$a$,$b$,$c$,$d$ are given by
$$
T^k_n\omega^{\alpha}_{\beta}=\omega^{\delta}_{\gamma}
\left(Q^k_m\right)^{\gamma\alpha}_{\beta\delta}T^m_n,\;\;\;
\mbox{where}\;\; Q=\lmatr{cc}A&B\\C&D\rmatr
\eqno(5.1)
$$
or
$$
a\omega^k=\omega^lA^k_la+\omega^lB^k_lc,~~~
c\omega^k=\omega^lD^k_lc+\omega^lC^k_la,
$$
$$
b\omega^k=\omega^lA^k_lb+\omega^lB^k_ld,~~~
d\omega^k=\omega^lD^k_ld+\omega^lC^k_lb,
$$

By using the commutation relations given above it is now easy to prove,
that
$$
(R_{12}T_1T_2-T_2T_1R_{12})\omega^k=0
$$
leads to the equations
$$
R_{12}Q_1Q_2=Q_2Q_1R_{12}
\eqno(5.2)
$$

It is an important property, that the quantum determinant be a central for
$\omega^k$ also
$$
\qdet\omega^k=(AD-qBC)\omega^k\qdet
\eqno(5.3)
$$
and we have
$$
AD-qBC=1
\eqno(5.4)
$$

Acting by the right exterior differential $\rdel$ on $(RTT-TTR)$ relation
and on the different expressions for the quantum determinant $\qdet$
$$
\rdel(ad-qbc) = (\omega^1 + \omega^kA^4_k - q\omega^kB^3_k)\qdet
$$
$$
\rdel(ad-qbc) = (\omega^1 + \omega^kA^4_k - q\omega^kC^2_k)\qdet
\eqno(5.5)
$$
$$
\rdel(da-q^{-1}bc) = (\omega^4 + \omega^kD^1_k - q{-1}\omega^kB^3_k)\qdet
$$
$$
\rdel(da-q^{-1}bc) = (\omega^4 + \omega^kD^1_k - q{-1}\omega^kC^2_k)\qdet
$$
we obtain the following consistency conditions
$$
C^2_k=B^3_k
$$
$$
\omega^3+\omega^kC^1_k-q^{-1}\omega^kA^3_k=0
\eqno(5.6)
$$
$$
\omega^2+\omega^kB^4_k-q\omega^kD^2_k=0
$$

By analogy with previous section we choose the representation of Q, where
$$
B=C=0,~~~ AD=1
\eqno(5.7)
$$
In this case we obtain the following anzatz for the matrices $A$ and $D$
$$
A=\lmatr{cccc}
1-t+\frac tu&0&0&(1-t)\frac tu\\
0&q&0&0\\
0&0&q&0\\
(u-1)&0&0&t\rmatr,~~~
D=\lmatr{cccc}
u&0&0&t-1\\
0&1/q&0&0\\
0&0&q&0\\
(1-u)\frac ut&0&0&(1-u+\frac ut)\rmatr,
\eqno(5.8)
$$
where $t(q)$ and $u(q)$ are some functions of $q$.

Using these matrices and applying right exterior differential $\rdel$ to
(5.1) together with the Cartan--Maurer equations for the
right--invariant 1-forms $\omega^k$ \ $\rdel\omega=\omega^2$ we can obtain
the following algebra of $\omega$-forms
$$(\omega^1)^2=(\omega^2)^2=(\omega^3)^2=(\omega^4)^2=0
$$
$$
\omega^3\omega^2+q^{-2}\omega^2\omega^3=0
$$
$$
\omega^1\omega^4+\omega^4\omega^1=0
$$
$$
u\omega^1\omega^3+(1-t+\frac tu)\omega^3\omega^1 =
(t-1)(q^2\omega^3\omega^4-\omega^4\omega^3),
\eqno(5.9)
$$
$$
(u-1)(q^{-2}\omega^1\omega^3-\omega^3\omega^1)=
q^2u\omega^3\omega^4+q^{-2}(1-t+\frac tu)\omega^4\omega^3,
$$
$$
u\omega^2\omega^1+(1-t+\frac tu)\omega^1\omega^2 =
(t-1)(q^2\omega^4\omega^2-\omega^2\omega^4),
$$
$$
(u-1)(q^{-2}\omega^2\omega^1-\omega^1\omega^2) =
q^2u\omega^4\omega^2+q^{-2}(1-t+\frac tu)\omega^2\omega^4
$$
and the condition
$$
t=u\;q^2
\eqno(5.10)
$$

Hence, taking into account (5.8), we get for  (5.5)
$$
\rdel\qdet= t(q^{-2}\omega^1+\omega^4)\qdet
\eqno(5.11)
$$

The parameter $t$ holds fixed form for the $\qtr\omega$.

By definition, the effect of applying right differential to an arbitrary
function $f$ on the quantum group is
$$
\rdel f=\rdel T^i_k\frac{\rpar}{\partial T^i_k}f=\omega^i\rnabla_if
\eqno(5.12)
$$
where $\rnabla_i$ are the right vector fields on the quantum group.

{}From the condition
$$
\rdel^2f=0
\eqno(5.13)
$$
and taking into account (5.9) we find a right vector field algebra
of arbitrary $t$.
$$
q^{-2}\rnabla_2\hrnabla_1-q^2\hrnabla_1\rnabla_2=(1+q^{-2})\rnabla_2,
$$
\nopagebreak
$$
q^{-2}\hrnabla_1\rnabla_3-q^2\rnabla_3\hrnabla_1=(1+q^{-2})\rnabla_3,
\eqno(5.14)
$$
\nopagebreak
$$
\rnabla_3\rnabla_2-q^{-2}\rnabla_2\rnabla_3=\hrnabla_1,
$$
$$
[\rnabla_4,\rnabla_2]= q^2(t-1)(q^2+1)\hrnabla_1\rnabla_2+
q^2(t+q^{-2}t-1)\rnabla_2,
\eqno(5.15)
$$
$$
[\rnabla_3,\rnabla_4]= q^2(t-1)(q^2+1)\rnabla_3\hrnabla_1+
q^2(t+q^{-2}t-1)\rnabla_3,
$$
where
$$
\hrnabla_1=\rnabla_1-q^{-2}\rnabla_4
\eqno(5.16)
$$

Again we have the decomposition of these commutation relations on the
$SL_q(2,C)$ and $U(1)$ subalgebras only for
$$
t=\frac{2}{1+q^{-2}},\;\;\;
u=\frac{2}{1+q^2}
\eqno(5.17)
$$

Commutation relations (5.14) form the algebra for $SL_q(2,C)$ and
for $U(1)$ we have
$$
[\hrnabla_4,\hrnabla_1]=[\hrnabla_4,\rnabla_2]=[\hrnabla_4,\rnabla_3]=0
\eqno(5.18)
$$
where $\hrnabla_4=\rnabla_1+\rnabla_4$.

After the mapping from $\hrnabla_1$,$\rnabla_2$,$\rnabla_3$,$\hrnabla_4$
to new generators $H$,$T_2$,$T_3$,$N$
$$
\hrnabla_1=\frac{1-q^{2H}}{1-q^2},\;\;
\rnabla_2=q^{H/2}T_2,\;\;
\rnabla_3=q^{H/2}T_3,\;\;
\hrnabla_4=N
$$
we again obtain the algebra
$U_qGL(2,C)$ in the form of the Drinfeld--Jimbo algebra (3.34).

Writing $\omega^1$ and $\omega^4$ as linear combinations
$$
\bar{\omega}^1=\frac{2}{q+q^{-1}}(\frac 1q\omega^1 + q\omega^4)=\qtr\omega,
\eqno(5.19)
$$
$$
\bar{\omega}^4=\frac{1}{1+q^{-2}}(\omega^1-\omega^4)
$$
we obtain more simple commutation relations for right 1--forms $\omega$
$$
\begin{array}{c}
[a,\bom^1]=[d,\bom^1]=0\\
\begin{array}{lcl}
a\bom^4=q^2\bom^4a, &~~~~~~~~& d\bom^4=q^{-2}\bom^4d\\
a\omega^2=q\omega^2a, && d\omega^2=q^{-1}\omega^2d\\
a\omega^3=q\omega^3a, && d\omega^3=q^{-1}\omega^3d
\end{array}
\end{array}
\eqno(5.20)
$$
The other relations are found by making the interchanges $a\to b$, $d\to
c$.

The commutation relations between $\omega$ have the following form
$$
\bom^1\omega^2+\omega^2\bom^1=0,~~~
\bom^1\omega^3+\omega^3\bom^1=0,~~~
\bom^1\bom^4+\bom^4\bom^1=0
\eqno(5.21)
$$
\nopagebreak
$$
q^2\omega^3\bom^4+q^{-2}\bom^4\omega^3=0,~~~
q^{-2}\omega^2\bom^4+q^2\bom^4\omega^2=0
$$
Thus we recovered the results [10] for right 1--forms.

Note, that commutation relations for right 1--forms $\omega$ can be
obtained from commutation relations for left 1--forms $\theta$ by making
interchanges ($a\leftrightarrow d$,$q\leftrightarrow 1/q$,
$\theta\leftrightarrow\omega$).

In the terms of the right differentials commutation relations between the
parameters and its differential have the following form for the
$GL_q(2,C)$
$$
\begin{array}{lcl}
a\;\rdel a=q^2\;\rdel a\;a+\frac{1-q^2}{2}\qtr\omega\;a^2\\
b\;\rdel b=q^2\;\rdel b\;b+\frac{1-q^2}{2}\qtr\omega\;b^2\\
c\;\rdel c=q^{-2}\;\rdel b\;b+\frac{1-q^{-2}}{2}\qtr\omega\;c^2\\
d\;\rdel a=q^{-2}\;\rdel d\;d+\frac{1-q^{-2}}{2}\qtr\omega\;d^2\\
b\;\rdel a=q\;\rdel a\;b+\frac{1-q^2}{2}\qtr\omega\;b\;a\\
a\;\rdel b=q\;\rdel b\;a+(q^2-1)\rdel
a\;b+\frac{1-q^2}{2}\qtr\omega\;a\;b\\
d\;\rdel c=\frac 1q \rdel c\;d + (q^{-2}-1)\rdel d\;c
+\frac{(1-q^{-2})}{2}\qtr\omega\;d\;c\\
c\;\rdel d=\frac 1q \rdel d\;c + \frac{(1-q^{-2})}{2}\qtr\omega\;c\;d\\
a\;\rdel c=q\rdel c\;a + \frac{(q-1/q)}{\qdet} (\rdel b\;c-q^{-1}\rdel
a\;d)\;a\;c + \frac{(1-q^{-2})}{2}\qtr\omega\;a\;c\\
c\;\rdel a=\frac 1q\rdel a\;c + \frac{(q-1/q)}{\qdet} (\rdel
b\;c-q^{-1}\rdel a\;d)\;c\;a + \frac{(1-q^{-2})}{2}\qtr\omega\;c\;a\\
a\;\rdel d=\rdel d\;a + \frac{(q-1/q)}{\qdet} (\rdel
b\;c-q^{-1}\rdel a\;d)\;a\;d + \frac{(1-q^{-2})}{2}\qtr\omega\;a\;d\\
d\;\rdel a=\rdel a\;d + \frac{(q-1/q)}{\qdet} (\rdel
b\;c-q^{-1}\rdel a\;d)\;d\;a + \frac{(1-q^{-2})}{2}\qtr\omega\;d\;a\\
b\;\rdel c=\rdel c\;b + \frac{(q-1/q)}{\qdet} (\rdel
b\;c-q^{-1}\rdel a\;d)\;b\;c + \frac{(1-q^{-2})}{2}\qtr\omega\;b\;c\\
c\;\rdel b=\rdel b\;c + \frac{(q-1/q)}{\qdet} (\rdel
b\;c-q^{-1}\rdel a\;d)\;c\;b + \frac{(1-q^{-2})}{2}\qtr\omega\;c\;b\\
b\;\rdel d=\rdel d\;b + \frac{(q-1/q)}{\qdet} (\rdel
b\;c-q^{-1}\rdel a\;d)\;b\;d + \frac{(1-q^{-2})}{2}\qtr\omega\;b\;d\\
d\;\rdel b=\rdel b\;d + \frac{(q-1/q)}{\qdet} (\rdel
b\;c-q^{-1}\rdel a\;d)\;d\;b + \frac{(1-q^{-2})}{2}\qtr\omega\;d\;b\\
\end{array}
\eqno(5.22)
$$

Again the choice $\qdet=1$,\ \  $\delta\qdet=\bom^1\qdet=0$
$$
\bom^1=\qtr\omega=\frac{2}{q+1/q}(\frac 1q\omega^1+q\omega^4)=0
\eqno(5.23)
$$
$$
\bom^4=\omega^1
$$
leads to the commutation relations for the $SL_q(2,C)$. It is not hard to
see that (5.35) with the conditions (5.36) have the form of the
solutions of Wess--Zumino for the quantum plane $C_q(2|0)$.

It will be noticed that now another combination of the parameters, namely
$$
ab=q\;ba,~~~~
cd=q\;dc
$$
and its differentials are quantum planes.

If we apply $\pi(b)=0$
$$
\begin{array}{lll}
T_-=\lmatr{cc}a&0\\c&d\rmatr=\lmatr{cc}y^{-1}&0\\x&y\rmatr\;\;\;
\end{array}
$$
and $\pi(c)=0$
$$
\begin{array}{lll}
T_+=\lmatr{cc}a&b\\0&d\rmatr=\lmatr{cc}x&y\\0&x^{-1}\rmatr\;\;\;
\end{array}
$$
so (5.22) corresponds to the solution of Wess--Zumino for the quantum
plane $C_q(2|0)$. Thus we see, that the right differential calculus
prefers the rows, but the left differential calculus prefers the columns.
At last we show the difference between matched differential calculus and
bicovariant differential calculus for $SL_q(2,C)$ group. Using commutation
relations (5.20) for right forms $\omega^k$, we can obtain commutation
relations for left forms $\theta^k$ after transformation $\theta =
g^{-1}\omega g$. For bicovariant differential calculus this commutation
relation must be the same as for left differential calculus. For matched
differential calculus we obtained commutation relations between the left
1--forms $\theta^k$ and higher degrees monomials, for example
$$
\theta^1b=\frac{1}{q^2}b\theta^1 - \frac{(q^4-1)}{q^4}b^2acd\theta^1-
\frac{(q^2-1)}{q^3}b^2ad^2\theta^3 +\frac{(q^2-1)}{q^4}b^2ac^2\theta^2
$$
$$
\theta^1c=q^2c\theta^1 + q^2(q^4-1)dc^2ba\theta^1+
q^2(q^2-1)dc^2b^2\theta^3 -q(q^2-1)dc^2a^2\theta^2
$$
$$
\theta^2b=\frac{1}{q^3}b\theta^2 - \frac{(q^4-1)}{q^5}b^3cd\theta^1-
\frac{(q^2-1)}{q^4}b^3d^2\theta^3 +\frac{(q^2-1)}{q^5}b^3c^2\theta^2
$$
$$
\theta^2c=qc\theta^2 + (q^4-1)cd^2ba\theta^1+
(q^2-1)cd^2b^2\theta^3 -\frac{(q^2-1)}{q}cd^2a^2\theta^2
$$
This solution is the solution of the equation (3.1) for matched
differential calculus also.

\section{Conclusions}

We proposed the construction of the differential calculus on the quantum
group and its subgroup with the property of the natural reduction: the
differential calculus on the quantum group  $GL_q(2,C)$ has to contain the
differential calculus on the quantum subgroup $SL_q(2,C)$ and
quantum plane $C_q(2|0)$ (''quantum matrjoshka''). We found, that there are
two differential calculi, associated to the left differential Maurer--Cartan
1-forms and to the right differential 1-forms. Matched reduction take the
degeneracy between the left and right differentials.

The space of 1-forms is four--dimensional for the quantum group
$GL_q(2,C)$ and is three--dimensional for the $SL_q(2,C)$. The quantum
determinant $\qdet$ is central element for 1-forms also. The obvious way
to carry differential calculus from $GL_q(2,C)$ over to $SL_q(2,C)$ by
imposing the determinant constraint $\qdet=1$ works with constraint
$\qtr\theta=\frac{2}{q+1/q}(q\theta^1+\frac 1q\theta^4)=0$.

Next step to carry differential calculus from $SL_q(2,C)$ over to $C_q(2|0)$
is to impose the constraint $\theta^2=0$ or $\theta^3=0$, that is
equivalent to $\pi(b)=0$ or $\pi(c)=0$.

The correspondence between left and right differential calculi is based on
well known property of $R$--matrices

$$
R_q = R^{-1}_{1/q}
$$
and by making the interchanges
$$
(a\to d,\;d\to a,\;q\to 1/q,\;b\to b,\;c\to c)
$$

The classical limit ($q\to 1$) of the ''left'' differential calculus and
the ''right'' differential calculus is the undeformed ordinary
differential calculus.

\section*{Acknoledgements}

We would like to thank professors D.Volkov, V.Drinfeld,J.Lukierski,
L.Vaksman, B.Zup\-nik, A.Isaev, P.Pjatov, A.Pashnev, A.Gumenchuk for
sti\-mu\-la\-ting dis\-cus\-sions and es\-pe\-ci\-al\-ly S.Krivonos for
the programm on the language ''Mathematica'' for the analitical
computation.One of the authors (V.G) would like to thanks Professor
J.Lukierski for the hospitality at the Institute of the Theoretical
Physics of the Wroclaw University,where part of this work was carried out.

This work was supported in part by International Science Foundation (grant
U21000) and Intermational Association for the promotion of Cooperation
with scientists from the independent states of the former Soviet Union
(grant INTAS-93-127).

\newpage

\section*{References}
\begin{enumerate}
\item V.G.Drinfeld {\sl Proc. Int. Congr. Math. (Berkeley, 1986)},
	p.798--810
\item Jimbo M. {\sf A $q$-difference analogue of $U_q(G)$ and the
	Yang--Baxter equation.}{\sl Lett.  Math.  Phys.} {\bf 10} (1985)
	p.63
\item Faddeev L., N.Reshetikhin, L.Takhtajan {\sl Algebra i
	Analiz},\mbox{\ }{ \bf\underline  1} (1989), p.178
\item Woronowicz S.L. {\sf Twisted $SU(2)$ Group. An Example of a
	non--commutative differential calculus.}{\sl Publ.  RIMS}, Kyoto
	Univ., vol.  {\bf 23}(1) 1987, p.117
\item Woronowicz S.L. {\sf Differential calculus on compact matrix
	pseu\-do\-groups\\ (quan\-tum groups).} {\sl Comm. Math.  Phys.}
	{\bf 122} (1989), p.125
\item J.Wess, B.Zumino {\sl Nucl.Phys.(Proc. Suppl.)}{\bf 18B}(1990),
	p.302
\item Yu.Manin {\sf Notes on quantum groups and quntum De Rham
	complexes}{\sl Teor.Mat.Fiz} v.{\bf 92} (1992) N3, p.425
\item Jur\v{c}o B. {\sf Differential calculus on quantized simple Lie
	groups.}{\sl Lett.  Math.  Phys.} {\bf 22} (1991), p.177
\item Faddeev L.D., Pyatov P.N. {\sf The Differential calculus on
Quantum Linear groups.} {\sl HEP-TH/9402070}

\item A.Schirrmacher, J.Wess, B.Zumino {\sl Z.Phys.} {\bf C49} (1990),
	p.317
\item Schupp P., Watts P. Zumino B. {\sf Differential geometry on linear
	quantum groups} {\sl Lett. Math. Phys. } {\bf 25} (1992), p.139
\item Sudbery A., \it Canonical differential calculus
        on quantum linear groups and supergroups,
        \rm Phys.Lett. \bf B284 \rm (1992)61-65
\item M\"uller--Hoissen F. {\sf Differential calculi on the quantum
	group $GL_{p,q}(2)$.} {\sl J. Phys. A: Math. Gen.} (1992), v.{\bf
	25}, p.1703
\item M\"uller--Hoissen F., Reuter C. {\sf Bicovariant differential
        calculi on $GL_{p,q}(2)$ and quantum subgroups.} {\sl J. Phys. A:
        Math. Gen.} (1993), v.{\bf 26}, p.2955

\item Isaev A.P., Pyatov P.N. {\sf $GL_q(N)$-covariant quantum algebras
	and covariant differential calculus.}{\sl Phys. Lett.}{\bf A189}
	(1993), p.81
\item Demidov E.E. {\sf Usp. Mat. Nauk} v.{\bf 48} (1993) N6, p.39--74
\item Ashieri P., Castellani L. {\sf An introduction to non--commutative
	differential geometry on quantum groups.} {\sl Preprint }
	CERN-6565/92
\item Zumino B. {\sf Differential calculus on quantum spaces and quantum
	groups} {Berkeley preprint LBL-33249 UCB-PTH-92141}
\item Schm\"udgen K., Schuller A.: {\sf Covariant differential calculi
	on quantum spaces and on quantum groups} {C. R. Acad. Sci.
	Paris}, v.{\bf 316} (1993), p.1155-1160
\item Sun X.D., Wang S.K. {\sf Bicovariant differential calculus on the
	two--parameter quantum group} {\sl J. Math. Phys.} {\bf 33}
	(1992), p.3313
\item Isaev A.P., Popowicz {\sl Phys. Lett.} {\bf B281} (1992), p.271
\item Akulov V.P., Gershun V.D., Gumenchuk A.I. {\sl Pis'ma Zh. Eksp.
	Teor. Fiz.}, {\bf 58}, No.6, p.462
\end{enumerate}

\end{document}